\documentclass[english]{article}
\usepackage{babel}
\usepackage{graphicx}
\title{A solution of Dirichlet problem using second partial derivatives of boundary function}
\author{Harry Yosh\\
APA Group, Australia\\
square17320508@yahoo.com}
\begin{document}
   \maketitle
\begin{abstract}
In Boundary Element Method, Green's function with no boundary conditions is used for solving Laplace's equation with Dirichlet boundary condition. To determine the gradient of solution on the boundary, we need to solve the boundary integral equation numerically in most practical cases.\\
Here we discuss the alternative method to avert solving that boundary integral equation. It is based on the solution for Poisson's equation which has the singularity on the boundary specified by the boundary function and shown it is applicable to arbitrary shaped domain and reduces calculation cost considerably.  \\ \\
{\bf keywords}\quad Boundary Element Method, Dirichlet problem, Poisson's equation, Green's function, boundary integral equation \\ \\
\end{abstract}

\paragraph{Introduction\\ \\}

In Boundary Element Method (BEM), Laplace's equation is solved by using the following Green's theorem;
\begin{eqnarray}
\int_D (\phi\nabla^2\psi - \psi\nabla^2\phi) dV = \int_{\partial D} (\phi\nabla\psi - \psi\nabla\phi) d\sigma
\end{eqnarray}
Let $\nabla^2\phi(x)=0$ and $\psi$ the Green's function $G$ for the Laplacian with no boundary conditions. Then,
\begin{eqnarray}
\phi(x) = \int_{\partial D} \left[G(x,x')\nabla'\phi(x') - \phi(x')\nabla' G(x,x')\right] d\sigma(x')
\end{eqnarray}
where $x$ is in the domain $D$. When $x$ is on the boundary $\partial D$ and $\partial D$ is sufficiently smooth around $x$, $\phi(x)$ satisfies the following equation;
\begin{eqnarray}
\frac{1}{2} \phi(x) = \int_{\partial D} \left[G(x,x')\nabla'\phi(x') - \phi(x')\nabla' G(x,x')\right] d\sigma(x') 
\end{eqnarray}
Under Dirichlet boundary condition, it is necessary to solve the Fredholm integral equation of the second kind shown in (3) for determining the gradient of solution on the boundary.$^{[1]}$

The boundary integral equation (3) is not solved analytically unless the boundary has some special forms$^{[4]}$ and usually solved numerically for practical purposes. It requires to deal with the large and fully populated matrix representing that equation for estimating the solution with sufficient precision, further the singularity caused from Green's function must be treated carefully.$^{[2]}$ Therefore generally the numerical solution of boundary integral equation requires great computational cost, and it is expected to improve the performance of BEM considerably by replacing this process with any alternative methods which don't require estimating the gradient of solution on the boundary. In the following section we discuss how to construct the method for satisfying such request. \\ \\

\paragraph{Second partial derivatives of boundary function for solving Dirichlet problem\\ \\} 

One of the methods to solve Dirichlet problem without estimating the gradient of solution on the boundary uses Poisson kernel. For example, the Dirichlet problem for the unit disk in the complex plane is solved by using Poisson kernel as,
\begin{eqnarray}
u(z)=\frac{1}{2\pi}\int_0^{2\pi} f(e^{i\theta}) \frac {1-\vert z \vert ^2}{\vert z-e^{i\theta}\vert ^2} d \theta \nonumber
\end{eqnarray} 
where $f(z)$ is the boundary function. Using this formula, the Dirichlet problems for various regions in the complex plane are solved by mapping unit disk to the region conformally$^{[3]}$ or applying the balayage method to the sequence of disks approximating the region.$^{[5]}$ However practically constructing such conformal map or applying the balayage method requires considerable effort when that region has complex shape. Here we discuss not only Poisson kernel, but also alternative boundary function in order to solve Dirichlet problem.

The above alternative boundary function is actually considered on the Poisson's equation $\nabla^2 \phi=g$ where the solution $\phi$ coincides with that of Laplace's equation $\nabla^2 \psi=0$ satisfying Dirichlet boundary condition: $\psi=f$ and $\phi=\psi$ in the domain surrounded by arbitrary shaped but sufficiently smooth boundary. Namely we are to look for the function $g$ which makes $\phi=\psi$ in that domain.

Firstly we discuss two-dimensional case in real plane ${\bf R}^2$. Let $f(x,y)$ the Dirichlet boundary condition of Laplace's equation $\nabla^2 \psi=0$, i.e. $\psi(x,y) = f(x,y)$ on the boundary $\partial D$ where $\partial D$ is sufficiently smooth. We suppose $D$ is closed and the boundary function is sufficiently smooth. For each point $(x,y)$ on the boundary $\partial D$ we take the local coordinate system $(x',y')$ where $x'$ is tangential on the boundary point and $y'$ is orthogonal to it. The origin of that coordinate system is the boundary point $(x,y)$. Also we suppose the Jacobian $\frac {\partial (x',y')}{\partial (x,y)}$ is equal with 1 (Fig.1). 

\begin{center}
  \includegraphics[height=60mm]{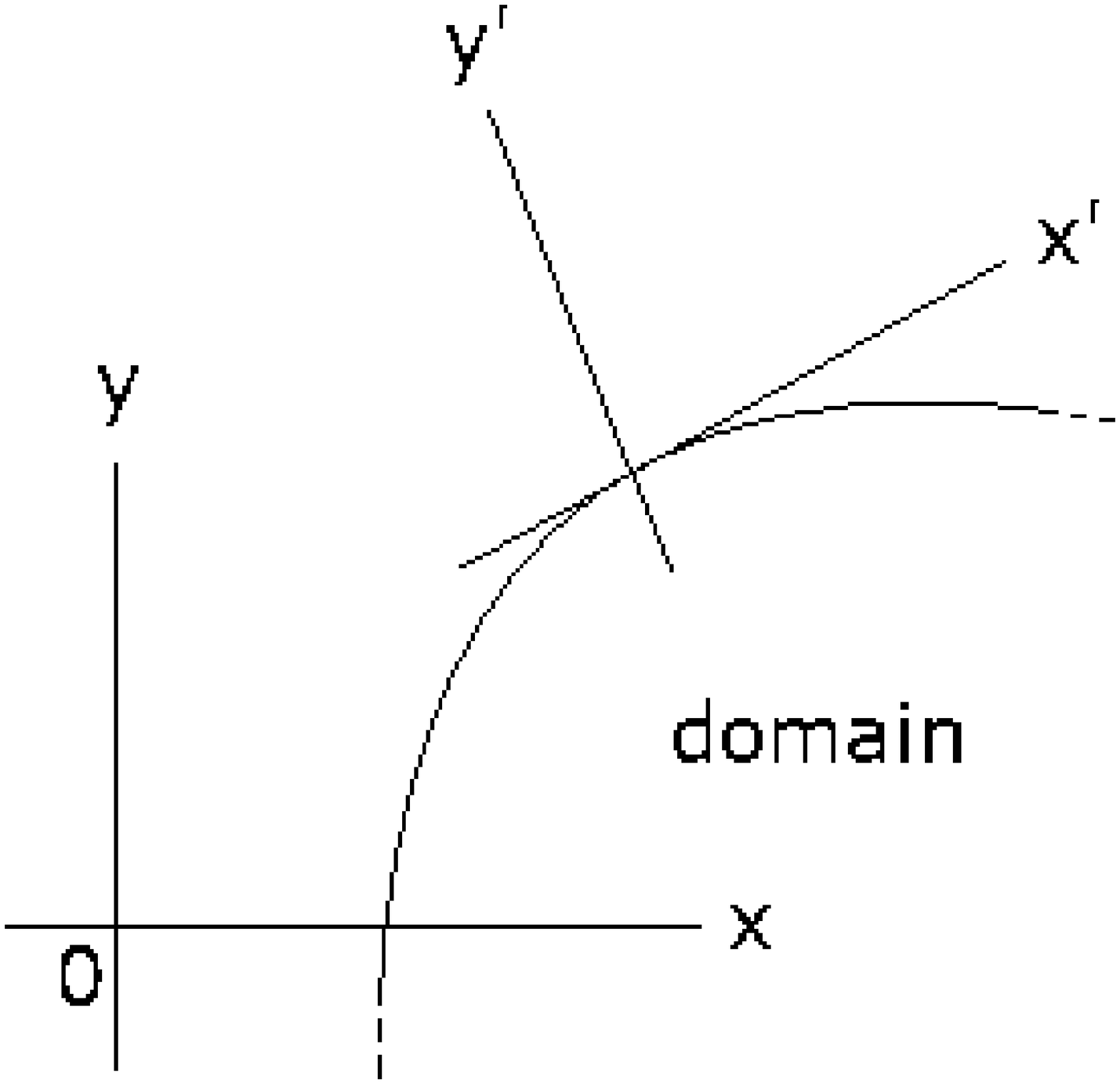}\\
  Fig.1 Local coordinate
\end{center}
$\hat f^{(2)}(x,y)$ is defined on the point ${\bf r}=(x,y) \in \partial D$ as,  
\begin{eqnarray}
\hat f^{(2)}(x,y)=\frac {\partial^2 f(x',y')}{\partial x'^2}|_{(x',y')=(0,0)}\quad ((x',y'):local\; coordinate)
\end{eqnarray} 
We suppose $\hat f^{(2)}(x,y)$ is also sufficiently smooth. Using $\hat f^{(2)}(x,y)$, $f^{(2)}(x,y)$ is defined as,
\begin{eqnarray}
f^{(2)}(x,y)=\hat f^{(2)}(x,y)\,\delta (y')
\end{eqnarray}
where $\delta (y')$ is one-dimensional delta function. Then the following proposition is satisfied; \\ \\
{\bf Proposition 1.}
\it Let $\psi$ the solution of Laplace's equation:
\begin{eqnarray}
\nabla^2 \psi=0
\end{eqnarray}
under the Dirichlet boundary condition $\psi(x,y) = f(x,y)$ on $\partial D$, then there exists $\phi$ such that,
\begin{eqnarray}
\psi = \phi + c \quad (c = constant) \quad : (x,y) \in D \nonumber
\end{eqnarray}
where $\phi$ satisfies the following Poisson's equation;
\begin{eqnarray}
\nabla^2 \phi(x,y)&=&f^{(2)}(x,y)\quad :(x,y) \in \partial D\; or \nonumber\\
&=&0 \quad \quad \quad \quad :(x,y) \not \in \partial D
\end{eqnarray} \\ \\ 
\rm
P r o o f. From the Poisson's equation (7), the Laplacian of $\phi$ at the point $(x,y)$ on the boundary is $f^{(2)}(x,y)$. Since $f^{(2)}(x,y)$ includes delta function, the Laplacian diverges on the boundary. However it does not mean that the second partial differential of $\phi$ along every direction always diverges on the boundary. To see it, firstly we define the function $g(x,y)$ as satisfying the following condition in the neighborhood of $(x,y)$.
\begin{eqnarray}
& &\nabla^2 g(x,y)=0\quad and \nonumber\\
& &\frac {\partial^2 g(x',y')}{\partial x'^2}|_{(x',y')=(0,0)}=\frac {\partial^2 f(x',y')}{\partial x'^2}|_{(x',y')=(0,0)}
\end{eqnarray}
where $(x',y')$ is the local coordinate around the point $(x,y)$, then a solution $\phi$ of Poisson's equation (7) is expressed in the neighborhood of $(x,y)$ as,
\begin{eqnarray}
\phi (x,y)&=&g(x',y') \quad \quad \quad \quad \quad \quad \quad \quad \quad \quad \quad \quad :if\; y'<0\; or \nonumber\\
&=&\frac {\partial^2 f(x',y')}{\partial x'^2}|_{(x',y')=(0,0)} \cdot y'+g(x',y') \;\; :if\; y' \geq 0
\end{eqnarray}
where the origin of the local coordinate is taken on $\partial D$. From the above expression, we know that the second partial differential of $\phi$ along $x'$ at the point $(x,y)$ on $\partial D$ is finite, i.e.
\begin{eqnarray}
\frac {\partial^2 \phi (x',y')}{\partial x'^2}|_{(x',y')=(0,0)}=\frac {\partial^2 f(x',y')}{\partial x'^2}|_{(x',y')=(0,0)}
\end{eqnarray}
Therefore $\phi$ is expressed in the neighborhood of $(x,y)$ fixing $y'=0$ as, 
\begin{eqnarray}
\phi(x',0)=f(x',0)+c_1x'+c_2\quad (c_1,c_2: constants)
\end{eqnarray}
Owing to the periodicity of $\phi(x,y)$ and $f(x,y)$ on $\partial D$, the constant $c_1$ in the above equation vanishes. On the other hand the solution $\breve \psi$ of Laplace's equation:
\begin{eqnarray}
\nabla^2 \breve \psi=0 \nonumber
\end{eqnarray}
under the Dirichlet boundary condition: $\breve \psi(x,y) = f(x,y)+c_2$ is expressed with $\psi$ in $(6)$ as,
\begin{eqnarray}
\breve \psi=\psi+c_2
\end{eqnarray}
Owing to the uniqueness of the solution of Dirichlet problem, $\breve \psi$ coincides with $\phi$ in the domain $D$, i.e.
\begin{eqnarray}
\phi=\psi+c_2
\end{eqnarray}
in $D$. {\quad //}\\

Since $f^{(2)}({\bf r})$ is the product of a smooth function and one-dimensional delta function as defined at (5), the solution $\phi$ of the Poisson's equation (7) expressed with Green's function as,
\begin{eqnarray}
\phi ({\bf r}) &=& \int_D f^{(2)}({\bf r'})G({\bf r}, {\bf r'})\, d\sigma'  
\end{eqnarray}
is actually the integral along the boundary $\partial D$. Therefore it is converted to one-dimensional integral using $\hat f^{(2)}({\bf r'})$ in (4) as,
\begin{eqnarray}
\phi ({\bf r}) = C \int_{\partial D} \hat f^{(2)}({\bf r'}) G({\bf r}, {\bf r'})\, dl'  
\end{eqnarray}
where $C$ is a constant (unit: length). If $\phi$ under the condition: $\phi({\bf r}) \rightarrow 0$ when ${\bf r} \rightarrow \infty$ is the solution, it must satisfy the following integral equation;
\begin{eqnarray}
& &\frac {\partial^2 \phi (x',y')}{\partial x'^2}|_{(x',y')=(0,0)} \nonumber \\
&=&\hat f^{(2)}({\bf r})=-\frac {C}{2\pi} \int_{\partial D} \hat f^{(2)}({\bf r'}) \frac {\partial^2 ln(|{\bf r'}-{\bf r}|)}{\partial x'^2}|_{(x',y')=(0,0)}\, dl'  
\end{eqnarray}
where ${\bf r}$ is taken on $\partial D$. Namely $-\frac {C}{2\pi}$ is the eigenvalue of the above integral equation.  As discussed later, $C$ is estimated numerically without solving that eigenvalue problem in the application to BEM.

The second partial differential of the boundary function $f(x,y)$ shown in (4) is actually carried out in the tangent space $T_p (\partial D)$ of each point on the boundary $\partial D$. Therefore when the boundary is not smooth enough to define the tangent space, the above proposition is not satisfied.\\

The above proposition is extended to higher dimensional one with the manner similar to two-dimensional case. In the n-dimensional case we suppose the boundary $\partial D$ is compact. The local coordinate $(x_1',...,x_{n-1}',x_n')$ at the point $p$ on the boundary is set so as to make $x_n'$ orthogonal to the tangent space $T_p(\partial D)$ and $\frac {\partial (x_1',...,x_n')}{\partial (x_1,...,x_n)}=1$. Also its origin is $p$. $\hat f^{(2)}({\bf r})$ is defined on the point ${\bf r}=(x_1,...,x_n) \in \partial D$ as,  
\begin{eqnarray}
\hat f^{(2)}({\bf r})=\sum_{j=1}^{n-1} \frac {\partial^2 f(x_1',...,x_n')}{\partial x_j'^2}|_{(x_1',...,x_n')=(0,...,0)}
\end{eqnarray} 
and $f^{(2)}({\bf r})$ is defined as,
\begin{eqnarray}
f^{(2)}({\bf r})=\hat f^{(2)}({\bf r})\, \delta(x_n')
\end{eqnarray} 
where $\delta (x_n')$ is one-dimensional delta function. Then the following proposition is satisfied;\\ \\
{\bf Proposition 2.}
\it Let $\psi$ the solution of Laplace's equation:
\begin{eqnarray}
\nabla^2 \psi=0
\end{eqnarray}
under the Dirichlet boundary condition $\psi({\bf r}) = f({\bf r})$ on $\partial D$, then there exists $\phi$ such that,
\begin{eqnarray}
\psi = \phi + c \quad (c = constant) \quad : (x,y) \in D
\end{eqnarray}
where $\phi$ satisfies the following Poisson's equation;
\begin{eqnarray}
\nabla^2 \phi({\bf r})&=&f^{(2)}({\bf r}) \quad :{\bf r} \in \partial D\; or \nonumber\\
&=&0 \quad \quad \quad :{\bf r} \not \in \partial D
\end{eqnarray} \\ \\
\rm
P r o o f. The function $g({\bf r})$ is defined as satisfying the following condition in the neighborhood of the point ${\bf r}=(x_1,...,x_n)$ in ${\bf R}^n$;
\begin{eqnarray}
& &\nabla^2 g({\bf r})=0\quad and \nonumber\\
& &\sum_{j=1}^{n-1} \frac {\partial^2 g(x_1',...,x_n')}{\partial x_j'^2}|_{(x_1',...,x_n')=(0,...,0)} \nonumber\\
&=&\sum_{j=1}^{n-1} \frac {\partial^2 f(x_1',...,x_n')}{\partial x_j'^2}|_{(x_1',...,x_n')=(0,...,0)}
\end{eqnarray}
Then a solution $\phi$ of Poisson's equation (21) is expressed in the neighborhood of $(x_1,...,x_n)$ as,
\begin{eqnarray}
\phi ({\bf r})&=&g(x_1',...,x_n') \quad \quad \quad \quad \quad \quad \quad \quad :if\; x_n'<0\; or \nonumber\\
&=&\sum_{j=1}^{n-1} \frac {\partial^2 f(x_1',...,x_n')}{\partial x_j'^2}|_{(x_1',...,x_n')=(0,...,0)} \cdot x_n' \nonumber\\
&+&g(x_1',...,x_n') \quad \quad \quad \quad \quad \quad \quad \quad \quad :if\; x_n' \geq 0
\end{eqnarray}
where the origin of the local coordinate is taken on $\partial D$. It satisfies the following equation;
\begin{eqnarray}
\sum_{j=1}^{n-1} \frac {\partial^2 \phi(x_1',...,x_n')}{\partial x_j'^2}|_{(x_1',...,x_n')=(0,...,0)}=\sum_{j=1}^{n-1} \frac {\partial^2 f(x_1',...,x_n')}{\partial x_j'^2}|_{(x_1',...,x_n')=(0,...,0)}
\end{eqnarray}
at the points on the boundary $\partial D$. Set $\sigma({\bf r})=\phi({\bf r})-f({\bf r})$, then the equation (24) is rewritten as a locally expressed Laplace's equation, i.e.
\begin{eqnarray}
& &\nabla^2  \sigma({\bf r}) \nonumber \\
&=&\frac {\partial^2 \sigma(x_1',...,x_n')}{\partial x_1'^2}+...+\frac {\partial^2 \sigma(x_1',...,x_n')}{\partial x_{n-1}'^2} \nonumber \\
&=&0
\end{eqnarray}
in the tangent space of each point ${\bf r}$ on the boundary $\partial D$.

As long as $\sigma({\bf r})$ is sufficiently smooth, $\sigma({\bf r})$ does not take maximum (or minimum) value in $\partial D$, otherwise at the maximum (or minimum) point the Laplacian $\nabla^2  \sigma({\bf r})$ would not be equal with zero. Since $\partial D$ is compact, $\sigma({\bf r})$ is a constant on $\partial D$, i.e.
\begin{eqnarray}
\phi({\bf r})-f({\bf r})=c \quad (c:constant)
\end{eqnarray}
On the other hand the solution $\breve \psi$ of Laplace's equation:
\begin{eqnarray}
\nabla^2 \breve \psi=0 \nonumber
\end{eqnarray}
under the Dirichlet boundary condition $\breve \psi({\bf r}) = f({\bf r})+c$ is expressed with $\psi$ in (20) as,
\begin{eqnarray}
\breve \psi=\psi+c
\end{eqnarray}
Owing to the uniqueness of the solution of Dirichlet problem, $\breve \psi$ coincides with $\phi$ in the domain $D$, i.e.
\begin{eqnarray}
\phi=\psi+c
\end{eqnarray}
in $D$. {\quad //}\\

If the boundary $\partial D$ is not compact, the solution of Laplace's equation (25) may take the maximum (or minimum) value on the closure of $\partial D$. Therefore it may not be constant and the proposition 2 is not satisfied in general. 

The solution $\phi$ of the Poisson's equation (21) expressed with n-dimensional integral as,
\begin{eqnarray}
\phi ({\bf r}) = \int_D f^{(2)}({\bf r'})G({\bf r}, {\bf r'})\, dV'  
\end{eqnarray}
is converted to (n-1)-dimensional integral as well as the two-dimensional case discussed in the proposition 1 as,
\begin{eqnarray}
\phi ({\bf r}) = C \int_{\partial D} \hat f^{(2)}({\bf r'})G({\bf r}, {\bf r'})\, d\sigma'  
\end{eqnarray}
where $C$ is a constant (unit: length). Actually $C$ is the eigenvalue of the following integral equation;
\begin{eqnarray}
& &\sum_{j=1}^{n-1} \frac {\partial^2 \phi (x_1',...,x_n')}{\partial x_j'^2}|_{(x_1',...,x_n')=(0,...,0)} \nonumber \\
&=&\hat f^{(2)}({\bf r})=C \sum_{j=1}^{n-1} \int_{\partial D} \hat f^{(2)}({\bf r'})\frac {\partial^2 G(x_1',...,x_n', {\bf r'})}{\partial x_j'^2}|_{(x_1',...,x_n')=(0,...,0)}\, d\sigma'
\end{eqnarray}
where ${\bf r}$ is taken on $\partial D$. \\ \\

\paragraph{Application to BEM\\ \\} 
 
From Proposition 2, the solution of Laplace's equation under Dirichlet boundary condition is expressed in general with the boundary integral (30) as,
\begin{eqnarray}
\psi ({\bf r}) = C_1 \int_{\partial D} \hat f^{(2)}G({\bf r}, {\bf r'})\, d\sigma' + C_2  
\end{eqnarray}
Since both ${\hat f}^{(2)}({\bf r})$ and $G(\bf {r}, \bf {r'})$ are explicitly calculable, the $\psi ({\bf r})$ expressed with (32) is estimated without any unknown factors except constants $C_1$ and $C_2$. Once we have got $C_1$ and $C_2$, we don't have to estimate the gradient of solution on the boundary, and it means the process for solving the boundary integral equation (3) has been eliminated. 

In order to apply the equation (32) to BEM, the boundary and boundary function must be finite and sufficiently smooth. The techniques for constructing mesh on the boundary used for ordinary BEM is effective as well. As the simplest way of numerical evaluation, the equation (32) is discretized as,
\begin{eqnarray}
{\tilde \psi}({\bf r}) = C_1\, \sum_{j=1}^{N} {{\tilde f}^{(2)}({\bf r}_j)\, G({\bf r}, {\bf r}_j)}\, \Delta \sigma({\bf r}_j) + C_2 
\end{eqnarray}
where ${\tilde f}^{(2)}$ is defined by discretizing (17) as,
\begin{eqnarray}
{\tilde f}^{(2)}({\bf r})=\frac {\Delta^2 f(x_1',...,x_n')}{\Delta x_1'^2}+\frac {\Delta^2 f(x_1',...,x_n')}{\Delta x_2'^2}+...+\frac {\Delta^2 f(x_1',...,x_n')}{\Delta x_{n-1}'^2}
\end{eqnarray} 
$C_1$ and $C_2$ are obtained by solving the following linear system;
\begin{eqnarray}
{\tilde \psi}({\bf r}_a)&=&C_1\, \sum_{j=1}^{N} {{\tilde f}^{(2)}({\bf r}_j)\, G({\bf r}_a, {\bf r}_j)}\, \Delta \sigma({\bf r}_j) + C_2 \\
{\tilde \psi}({\bf r}_b)&=&C_1\, \sum_{j=1}^{N} {{\tilde f}^{(2)}({\bf r}_j)\, G({\bf r}_b, {\bf r}_j)}\, \Delta \sigma({\bf r}_j) + C_2 
\end{eqnarray}
where ${\bf r}_a$ and ${\bf r}_b$ are the mesh points chosen arbitrarily on the boundary. ${\tilde \psi}({\bf r}_a)$ and ${\tilde \psi}({\bf r}_b)$ are given previously as the boundary condition. 

Since the calculation loads for obtaining each ${\tilde f}^{(2)}({\bf r}_j)$, $G({\bf r}, {\bf r}_j)$, and $\Delta \sigma({\bf r}_j)$ don't depend on the number of mesh points $N$, the order of calculation load for estimating $\psi({\bf r})$ inside of the boundary is $O(N)$.

In the Fig.2, two-dimensional Laplace's equation solved by the above method is illustrated. The Laplace's equation dealt with here has the following boundary conditions;
\begin{eqnarray}
& & \partial D = \{(x,y)|x^2+y^2=1\} \nonumber\\
& & f(cos\, \theta, sin\, \theta)=sin\, 2 \theta \quad (0 \leq \theta < 2\pi) \nonumber
\end{eqnarray} 

\begin{center}
  \includegraphics[height=90mm]{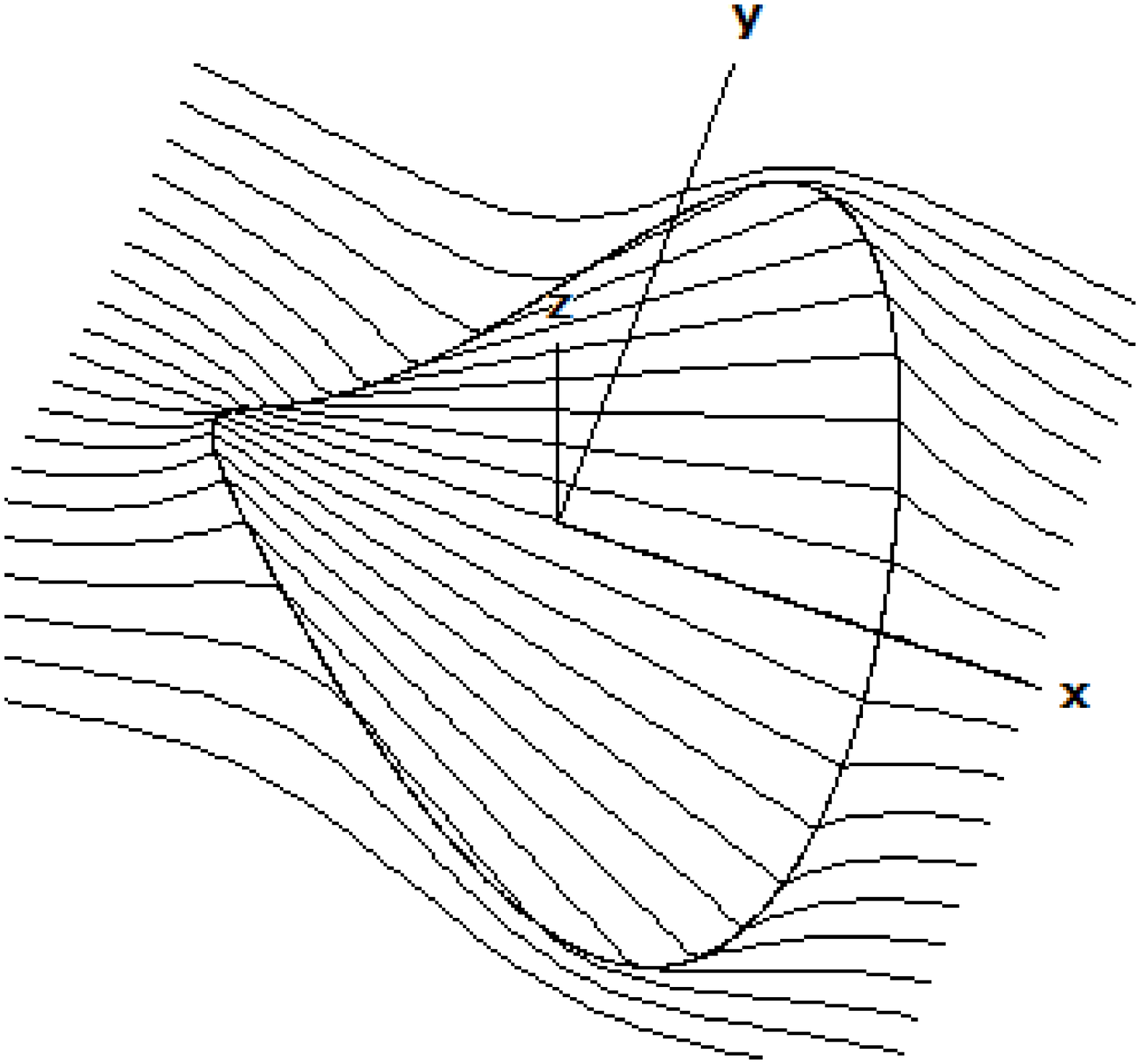}\\
  Fig.2 Numerical evaluation of Poisson's equation having singularity along the boundary
\end{center}
Actually it is the numerical solution of Poisson's equation (7) and defined not only in the domain $D$ for the Laplace's equation but also outside of it. As seen in the figure, the second partial differentials of the solution diverge on the boundary except the direction tangential to the boundary as mentioned in the proof of Proposition 1.\\ \\

\paragraph{Summary\\ \\} 
  
The solution of Laplace's equation under Dirichlet boundary condition coincides with the solution of a Poisson's equation which has a specific singularity along the boundary. The domain in which the solution of that Poisson's equation is applicable is arbitrary as far as the boundary of the domain is compact and the boundary function is finite and both of them are sufficiently smooth. 

Comparing the solution derived from the above Poisson's equation with that from ordinary Boundary Element Method, we find the former does not need to solve the boundary integral equation which is used for estimating the gradient of the solution on the boundary and requires considerable calculation cost in the latter. Therefore it is expected to reduce the calculation cost by implementing the method discussed here to Boundary Element Method especially when dealing with the Laplace's equation in the domain having complex shape. \\ \\


\begin{thebibliography}{5}
\bibitem{Wrobel}
      L. C. Wrobel, M. H. Aliabadi, 
      \emph{The Boundary Element Method}.
      New Jersey: Wiley, 2002.
\bibitem{Beer}
      G. Beer, I. Smith, C. Duenser,
      \emph{The Boundary Element Method with Programming: For Engineers and Scientists}.
      Springer-Verlag, 2008.
\bibitem{Conway}
      J. B. Conway,
      \emph{Functions of One Complex Variable I}.
      Springer-Verlag, 1978.
\bibitem{Axler}
      S. Axler, W. Ramey, 
      \emph{Harmonic polynomials and Dirichlet-type problems}.
      Proc. Amer. Math. Soc. 123. 3765-3773, 1995.
\bibitem{Poincare}
      H. Poincar$\acute e$, 
      \emph{Sur les $\acute e$quations aux d$\acute e$rivees partielles de la physique math$\acute e$matique}.
      Amer. J. Math., 12 : 3, pp. 211-294, 1890.
\end{thebibliography}
\end{document}